\documentclass[10pt]{article}

\usepackage[]{amsmath,amscd,amsfonts,amsthm,amssymb}
\usepackage[greek, english]{babel} 
\languageattribute{greek}{polutoniko}

\usepackage{graphicx}
\usepackage[all]{xy}
\usepackage{psfrag}
\usepackage[lflt]{floatflt}
\newcommand{\bq}{\begin{quote}}
\newcommand{\eq}{\end{quote}}

\newtheorem{Th}{Theorem}
\newtheorem{ax}{Axiom}
\newtheorem{lm}{Lemma}
\newtheorem{df}{Definition}
\newtheorem{pr}{Proposition}
\newtheorem{cl}{Corollary}
\newtheorem{re}{Remark}
\newtheorem{as}{Assumption}
\newtheorem{wg}{Wild Guess}
\newtheorem{ex}{Example}
\newcommand{\bth}{\begin{Th}\hspace{-5pt}{\bf .} \ }
\newcommand{\Eth}{\end{Th}}
\newcommand{\bax}{\begin{ax}\hspace{-5pt}{\bf .} \ }
\newcommand{\eax}{\end{ax}}
\newcommand{\blm}{\begin{lm}\hspace{-5pt}{\bf .} \ }
\newcommand{\elm}{\end{lm}}
\newcommand{\bdf}{\begin{df}\hspace{-5pt}{\bf .} \ }
\newcommand{\edf}{\end{df}}
\newcommand{\bpr}{\begin{pr}\hspace{-5pt}{\bf .} \ }
\newcommand{\epr}{\end{pr}}
\newcommand{\bcl}{\begin{cl}\hspace{-5pt}{\bf .} \ }
\newcommand{\ecl}{\end{cl}}
\newcommand{\bre}{\begin{re}\hspace{-5pt}{\bf .} \ }
\newcommand{\ere}{\end{re}}
\newcommand{\bas}{\begin{as}\hspace{-5pt}{\bf .} \ }
\newcommand{\eas}{\end{as}}
\newcommand{\bwg}{\begin{wg}\hspace{-5pt}{\bf .} \ }
\newcommand{\ewg}{\end{wg}}
\newcommand{\bex}{\begin{ex}\hspace{-5pt}{\bf .} \ }
\newcommand{\eex}{\end{ex}}

\newcommand{\bit}{\begin{itemize}}
\newcommand{\eit}{\end{itemize}\par\noindent}
\newcommand{\ben}{\begin{enumerate}}
\newcommand{\een}{\end{enumerate}\par\noindent}
\newcommand{\beq}{\begin{equation}}
\newcommand{\eeq}{\end{equation}}
\newcommand{\beqa}{\begin{eqnarray*}}
\newcommand{\eeqa}{\end{eqnarray*}\par\noindent}
\newcommand{\beqn}{\begin{eqnarray}}
\newcommand{\eeqn}{\end{eqnarray}\par\noindent}
\usepackage{fancybox}

\title{\huge Infinity and the Sublime\large$^{\footnote{This paper is dedicated to the memory of my father, Emiel Verelst ($^\circ$26/12/1934 -- $\dagger$21/10/2009).}}$} 
\author{} 
\date{}

\begin{document}   

\maketitle  

\vspace{-1cm} 

%\centerline{\small \sc {[DRAFT]}}

\centerline{Karin Verelst}
\vspace{1mm}
\par
{\scriptsize
%\centerline{{\em FUND-CLEA}}
\vspace{1mm}
\centerline{{\em Vrije Universiteit Brussel\footnote{FUND-CLEA, Dept. of Mathematics.} and}}
\centerline{{\em Open Universiteit Nederland\footnote{Faculty of Cultural Studies.}}}
\centerline{{\em kverelst@vub.ac.be}}
} 

\bigskip
\bigskip
\bigskip

\begin{flushright}{\footnotesize {\em Nommer, c'est avoir individu}\\{Nikolai N. Luzin}}
\end{flushright}
\bigskip

\subsection*{\sc Abstract}
%\bigskip
\frenchspacing In this paper we intend to connect two different strands of research concerning the origin of what I shall loosely call ``formal'' ideas: firstly, the relation between logic and rhetoric - the theme of the 2006 Cambridge conference to which this paper was a contribution -, and secondly, the impact of religious convictions on the formation of certain twentieth century mathematical concepts, as brought to the attention recently by the work of L. Graham and J.-M. Kantor. In fact, we shall show that the latter question is a special case of the former, and that investigation of the larger question adds to our understanding of the smaller one. Our approach will be primarily historical.
\subsection*{\sc Introduction}
%\bigskip

\noindent  Let us start with an outline of the issues that are at stake. In the conclusion to their {\em Isis}-paper, Graham and Kantor state: 

 \begin{quote} {\em We believe that our study of French and Russian developments in set theory and the theory of functions points strongly toward the importance of cultural factors in the process and creation of mathematics --- in the French case, Descartes, positivism, and Pascal; in the Russian case, mystical religious beliefs, particularly those of the Name Worshipping movement. As a result the French and the Russians followed different approaches.\footnote{L. Graham, J.-M. Kantor, ``A comparison of two cultural approaches to mathematics France and Russia, 1890-1930'',  {\em Isis,} {\bf 97}, 1, 2006.} }\end{quote}

\noindent They immediately add that such ``intellectual causation'' can never be proved rigorously, only made plausible, and they are evidently right. But it is equally evident that this holds true for many phaenomena investigated in, loosely speaking again, cultural studies, without making them any less relevant. Graham and Kantor also point out that comparative studies add up to more credible evidence, and to this we also agree. Therefore we shall in the present paper enlarge the comparative scope of their work into the past, and see whether the conclusions they reached do still stand. We limit ourselves to the ``French'' part of their comparative equation, not only because we possess the linguistic and historical tools required to do the concerned research properly, but also because our authors elaborate the French background only little when compared to the Russian one, probably because they assume it to be ``known''. However, when looked into more deeply, some surprises are bound to turn up! The starting point for our investigation is fixed by our authors when they say that the French intellectual framework with respect to ``rationalism'' has been shaped decisively by Descartes and Pascal.  This brings us to the beginnings of the Early Modern period, a period of cultural, religious and scientific turmoil indeed. It is our intention to show in the following pages that during that period a debate --- really more like an intellectual war --- was going on concerning the nature of rhetoric, and especially memory, in which the different positions involved had fundamentally different theological roots, and lead to different practices in philosophy, art (in its present-day sense), and science, which are relevant to the Graham-Kantor thesis. One of the focal points of the debate was the notion of {\bf \em method}. The final victory of one of the intellectual parties in the debate brought about the destruction of the other methodological approach and stood godmother to the birth of modern science, in its rationalist outlook advocated by Descartes, Leibniz and Pascal, in its empirist cloak by Newton and his followers. It is relevant for our concerns that the main characters in the final stage of this heroic tale gave birth to the kind of mathematics repudiated by the masters of the later Russian mathematical school, viz., infinitesimal calculus. 

With respect to ``method'', it is a commonplace in our philosophical tradition that Descartes initiated a new way of doing ``natural philosophy'' --- for him unquestionably aequivalent to doing science ---, by the introduction of a new method based on ``clear ideas'' and mathematical reasoning. The title of his famous book, {\em Discours de la m\'{e}thode}, brings the point home explicitly. This literally methodological emphasis framed his approach to natural philosophy from the start, as is clear from the title of his first, but posthumously published, work {\em Regulae at directionem ingenii} ({\em Rules for the Direction of the Mind}).\footnote{R. Descartes, {\em Regulae ad directionem ingenii}, Texte critique, G. Crapulli (ed.), Martinus Nijhoff, La Haye, 1966. (This edition is better than the version in Adam et Tannery (vol. X, p. 362 sq., Cerf, Paris, 1897-1913), because it takes into account the Dutch translations published after 1680.} This approach became the intellectual program of the period, fostered not only by Descartes, but also by the famous {\em Port Royal Logic}. Arnauld and Nicole, authors of the {\em Port Royal Logic}, call it {\em L'Art de Penser}.\footnote{A. Arnauld et P. Nicole, {\em La Logique ou L'Art de Penser}, Flammarion, Paris, 1970, [1662].} This very influential and widespread book outlines an anti-Aristotelian conception of logic based on the study of natural language and on the nature of mathematical deductions (as opposed to merely syllogistic ones), and it states explicitly that the authors chose their examples in such a way as to introduce the reader to the correct philosophy of nature, as well as to correct methods for judgment in matters moral and religious. Arnauld and Nicole heavely stress the importance of both simplicity and 'good sense' in judgment, and lay this down in eight `Rules' destined to govern the process of reasoning. These are exemplified in the `Quatri\`{e}me partie. De la m\'{e}thode', the last part of their book. They furthermore distinguish carefully between {\em analyse}, the method of invention, and {\em synthese}, the method of composition.\footnote{A. Arnauld et P. Nicole, {\em o.c.}, p. 368.} Finally, they declare explicitly that in this respect they are followers of Descartes in his {\em Discours de la m\'{e}thode}. The example was followed by Newton\footnote{Newton had a copy of the latin translation of the Port Royal Logic in his personal library: A. Arnauld \& P. Nicole, {\em Logica, sive Ars cogitandi: in quae praeter vulgares regulas plura nova habentur ad rationem dirigendam utilia}, E. Tertia apud Gallos ed. recognita \& aucta in Latinum versa, Londini, 1687.  In the Harrison catalog of Newton's library, this is item 79. I have reasons to believe it has been given or pointed out to him by his friend Nicolas Fatio de Duillier.}, who starts the third Book of his {\em Principia} --- not incidentially named after another work by Descartes, the {\em Principia Philosophiae}\footnote{A link evident to his contemporaries, as witnessed by Huygens in a letter to Leibniz dd. May 1694, where he refers to the ``raisonnement et experiences de Newton dans ses Principes de Philosophie'', {\em Oeuvres Compl\`{e}tes de Christian Huygens},  X, n¡ 2854, p. 614. On the relation between Newton and Descartes, see, e.g., A. Janiak, {\em Newton as Philosopher}, Cambridge University Press, Cambridge, 2008.} --- with a notorious set of (different!) ``Regulae philosophandi'', rules to which scientific practice should comply in order to lead to trustworthy new results.\footnote{The genesis of Newton's ``rules'' is problematic in itself and deserves  far more attention than it usually gets. A notable exception is I.B. Cohen, ``Hypotheses in Newton's Philosophy'', {\em Physis}, VIII, 1966. An interesting study of scope and aim of Newton's empiricism as voiced in their final version is: A. Shapiro, ``Newton's Experimental Philosophy'', {\em Early Modern Science and Medicine}, {\bf 9}, 3, 2004.} 

But Descartes had his immediate predecessors as well. Five years before the publication of the {\em Discours}, a small group of scholars gathered to discuss  the then current ideas on ``method'', and published the transactions of this private academy under the title ``De la m\'{e}thode''.\footnote{F.A. Yates, {\em The Art of Memory}, University of chicago Press, Chicago, 1966, p. 370, ft. 4.} Shortly before that, a notorious critic and reformer of the scholastic curriculum at the universities of the time, Peter Ramus, had used the term to refer to his new, ``dialectical'' system of memory, designed to replace the older ``rhetorical'' one. Yates points out that this in itself suggests that there might be a more intrinsic relation between ``the history of memory and the history of method''.\footnote{F.A. Yates, {\em o.c.}, p. 369.} In what follows we take this suggestion seriously and show that this supposedly original theme reaches back into the distant past and connects the birth of modern science in an unexpected way to the older philosophical tradition. The connection runs through the notion of {\em inventio}, the art of finding ``new ideas'', new and truthful insights on any subject of human interest to be presented before an audience or readership.\footnote{The difference was vague well into the XVI-th century.  During the Middle Ages, books were considered as living beings to almost the same extend as the human beings who carried them around physically or in their memory, ``borne everywhere in the minds of their {\em listeners}'', R. de Bury,  {\em Philobiblon}, 4 (my italics). For more on this, see below.} Both the origins and repercussions of the different options regarding {\em inventio} have much broader cultural and religious implications. We shall  furthermore see how, dually, different choices with respect to {\em inventio} appear to determine the whole intellectual outlook of the individual committed to them. Matters of `style' can act as light-house fires on conceptual vessels searching for the metaphysical home ports in which they were constructed, and by means of which they ought to be construed: {\em Les questions controvers\'{e}es de rh\'{e}torique (...) jouent donc un r\^{o}le d\'{e}terminant dans lÕEurope catholique. Elles t\'{e}moingent de pr\'{e}f\'{e}rences qui vont quelquefois jusqu'\`{a} des choix th\'{e}ologiques tranch\'{e}s.}\footnote{  M. Fumaroli, L'\'{E}cole du silence, Le sentiment des images au XVIIe si\`{e}cle, Flammarion, Paris, 1998 [1994], p. 16.} This holds  true not only for the catholic world, and remains so long after the advent of Modernity, so that issues related to {\em inventio} provide us with a reliable compass to follow the often hardly visible lines of fracture along which certain debates concerning the foundations of mathematics have been taking place at the verge of the twentieth century. We shall show how in origin rhetorical conceptions feeding on metaphysical and even theological grounds continue to inform present day mathematical thought, in particular in France, in a way that confirms the findings of Graham and Kantor, but sheds a different light on them. 

In order to set out the different positions in the debate, we focus on the notion of the {\em sublime}, a concept central to early modern controversies on the origin of rhetorical creativity. The sublime surpasses the rule-governed {\em inventio}, as if it were a divine breath that infuses inspiration into the mind, causing a {\em furor} or {\em enthusiasm} in the mind affected by it. According to Adrien Baillet, his biographer, Descartes was well aware of the reality and force of this imaginative enthusiasm:

\begin{quote} {\em [Descartes] ne croioit pas qu'on d\^{u}t s'\'{e}tonner si fort de voir que les Po\`{e}tes, m\^{e}me ceux qui ne font que niaiser, fussent pleins de sentences plus graves, plus sens\'{e}es, \& mieux exprim\'{e}es que celles qui se trouvent dans les \'{e}crits des Philosophes. Il attribuait cette merveille \`{a} la divinit\'{e} de l'Enthousiasme, \& \`{a} la force de l'imagination.\footnote{From {\em Vie de Monsieur Des-Cartes} (1691), quoted in: M.H. Keefer, ``the Dreamer's Path: Descartes and the Sixteenth Century'', {\em Renaissance Quarterly}, {\bf 49}; 1996, p. 30.}} \end{quote} 

\noindent This possibility of a direct contact between finite minds and the divine realm was prone to be considered a blasphemy, and it crumbled under the rationalising weight of the Reformation, both protestant and catholic, that sought to eradicate the last traces of ``paganism'' in all aspects of European cultural life: 

\begin{quote} {\em D\`{e}s lors que l'enthousiasme, soumis \`{a} la critique patiente et sys\-t\'{e}\-ma\-tique de l'humanisme \'{e}rudit, depuis Jules-C\'{e}sar Scaliger jusqu'\`{a} M\'{e}ric Casaubon, n'apparait plus comme le principe s\'{e}minal de la connaissance et de l'invention humaines, d\`{e}s lors que le principe de raison s'impose \`{a} sa place, ce sont des pans entiers de la culture humaniste qui s'\'{e}croulent, pour faire place \`{a} un nouvel \'{e}difice de style moderne.}\footnote{M. Fumaroli, ``Le cr\'{e}puscule de l'enthousias\-me'', in: {\em H\'{e}ros et Orateurs,} Droz, Gen\`{e}ve, 1996, p. 376.} \end{quote} 

\noindent  Our intention is to show that this decline ultimately had consequences for formal ``arts'' like mathematics as well. Acceptance or rejection of the sublime in rhetoric is connected to the belief in an absolute or infinite realm to which a great soul can connect, and has far-reaching theological consequences. It cuts through the frontlines of the notorious ``querelle des Anciens et Modernes'', and bears upon precise conceptions of the relation between language, existence and imagination. The sublime paves the way to the acceptance or rejection of the notion of mathematical infinity. This intellectual connection is straightforward and can be shown to shape. e.g.,  certain oppositions between Descartes and Pascal, but it remains operative well into the twentieth century, where it explains the different approaches developed to the actual infinite in the foundations of mathematics by the French rational - heir to Descartes and Port Royal - school (Poincar\'{e}, Lebesgue, Hadamard, Baire, Borel), and by the Moscow school of mathematics (Egorov, Luzin, Florenskii, Suslin), feeded by a strand of orthodox mysticism called the ``Name-worshippers'', who shared a central tenet of the defenders of the sublime, c.q., that , at least to a certain extend, `to name' and `to exist' coincide.\footnote{Apart from the paper already cited, consult L. Graham and J.-M. Kantor, ``Name Worshippers. Russian religious mystics and French rationalists: Mathematics 1900-1930'', {\em American Academy of Arts and Sciences Bulletin}, vol. LVIII, nr 3, 2005. See also their recent book: L. Graham and J.-M. Kantor, {\em Naming Infinity. A True Story of Religious Mysticism and Mathematical Creativity}, Belknap Press, Cambridge, Massachusetts, 2009.}

%\newpage
\subsection*{\sc The Rhetorical Structure of the World}

We are so used to think by ourselves that we almost forget how rare and strange this psychological feat actually is. Yet the ability to retract `within' oneself is a recent and culturally determined phaenomenon which was hardly available throughout the major part of human history. Indeed, to think was to think aloud. In our own Antiquity to think and to speak did largely coincide.\footnote{R.B. Onians, {\em The Origins of European Thought}, Cambridge University Press, Cambridge, 1994 [1951], p. 13}  If we look at one of the traditional definitions of rhetoric as the {\em ars bene dicendi} with this in mind, then the equally traditional viewpoint that moores rhetoric so firmly to the correct ways to think, i.e. the {\em ars inveniendi,} the method for discovery by which the mind ``finds'' new ideas, becomes something less of a surprise. This ``art of discovery'' was itself linked to the {\em ars memorativa}, the use of memory, that vast storehouse for both knowledge and imagination, the groundplan of which was organised by means of common places as well as particular ones, and seemed to have the awesome capacity of bringing order and structure into the world.\footnote{The pioneering work studying these interplays during Renaissance is the already cited book by F.A. Yates, {\em The Art of Memory}. More recent contributions to the field include L. Bolzoni, {\em La stanza della memoria}, Einaudi, Turino 1995 (translated into English as {\em The Gallery of Memory} in 2001, and M. Carruthers, {\em The Craft of Thought}, Cambridge University press, Cambridge, 1998 (for the Middle Ages).} Rhetoric implies much more than mere eloquence. This moreover works in both ways: in the Mediaeval period, books were seen as extensions of living memory, not the other way around: {\em (...) books demanded to be alive, to speak and converse, to be consumed and digested through the memories of living people.}\footnote{M. Carruthers, ``Mechanisms for the Transmission of Culture: The Role of ``Place'' in the Arts of Memory'',  in: {\em Translation, or the Transmission of Culture in the Middle Ages and the Renaissance}, L. H. Hollengreen (ed.), Brepols, Turnhout, 2008.} Up to early modernity it was common to label {\em Rhetorica}\/ treatises on eloquence, but also those on matters spiritual or metaphysical. Now a crucial rhetorical notion was that of {\em method}, the way to be followed by the mind in order to discover new ideas, as well as the way to present them appropriately. This was discussed in a work ascribed to the Hellenistic rhetor Hermogenes: {\em On the Method of Forcefulness}.\footnote{G.A. Kennedy (ed.), {\em Invention and Method: Two Rhetorical Treatises from the Hermogenic Corpus}, Brill, Leiden, 2005, p. 414.} Thus when Descartes, that fierce opponent to the rhetorical ways of thinking prevalent during his lifetime, nevertheless entitled one of his most famous books {\em Discours de la m\'{e}thode}, he stayed faithful to a vestige of that same tradition. Our present disdain for rhetoric is the exception, and one that becomes comprehensible only in the context of the history of rhetoric itself. The strong interconnections between the religious and cultural reformatory fervour that hit Renaissance Europe during the sixteenth century, and the emergence of new conceptions on the functioning of the human mind, have been brought to light in all their stunning complexity by F.A. Yates in her remarkable studies on Renaissance memory. But memory was only one of the five traditional canons of rhetorical art. This is why the Renaissance revolution initiated by Petrus Ramus with respect to the use of memory articulates with many other fundamental intellectual transitions during that period. Before turning our attention to memory, we have to revise shortly the major steps in the more general process.\\

During the Middle Ages, the rhetorical works by Cicero and Quintilian had been accesible only incompletely, while the Greek sources, except Aristotle, had been lost. Cicero's {\em De Inventione} and the {\em Ad Herrenium} (erroneously attributed him) were very influential treatises, in which the basic techniques for invention and rhetorical exposition were given in a textbook style. In academical circles, however, the Aristotelian tradition with respect to rhetorics as separated from and inferior to dialectics reigned supreme: the focus was on {\em truth} as the most important means of persuasion. In this the scholastics were followers of Boethius, according to whom dialectic comes down to logic: its aim is to ascertain the truth of a statement by controversial argument, organised through the use of {\em common places} (definitions, genera, species, qualities \&.), which therefore belong to the realm of dialectics alone. This view starts to change when in 1415, Poggio Bracciolini discovered in the monastery of St. Gall a complete copy of Quintilian's {\em De Institutione oratoria}.\footnote{T.M. Conley, {\em Rhetoric in the European Tradition},  University of Chicago Press, Chicago, 1990, p. 112 sq..} Its key notions, such as the core relation between virtue and eloquence and the {\bf universal range of rhetoric}, caused a stirr throughout Renaissance Italy. The discovery of the complete Ciceronian {\em De Oratore} by G. Landriani in Milan added to the commotion. Here we had a text written in an extremely beautiful style, by means of dialogues, not by means of dogmatic exposition, utterly different from the dry schemes of the School. This forms the backdrop to a conflictual development that in fact had started already with Petrarc, cut through the religious and cultural frontlines in Europe from tne XVthe to the XVIIIth century and which became known, particularly in France, as {\bf La Querelle des Ancies et Modernes}. Evidently, the debate shifted in shape, place and focus throughout that period: the scholastic ``moderns'' which were the opponents to Erasmus differed in many respects from the Cartesian rationalists who were the ennemies of Montaigne, La Fontaine and Mme Dacier, or the empirists who rejected Burke. Nevertheless, the line of fracture connecting the parties at both sides of the breach remained the same: {\em Les uns veulent arrimer l'Europe Moderne au g\'{e}nie antique. Les autres veulent s'en \'{e}manciper.}\footnote{M. Fumaroli, {\em La Querelle des Anciens et des Modernes}, Gallimard, Paris, 2001, p. 8.} Understandebly enough, in the Italian city states of the XIV-XVth centuries, the ideal of the orator-statesman had an evident appeal. Cicero thaught 1) that the primary duty of man is action, and 2) that the vita activa does not distract from one's intellectual powers, but on the contrary stimulates them. Humanists in, say, Florence were politically active men. Lorenzo Valla, former chancellor to the Republic, was the author of two major works on rhetoric: {\em Elegantiae linguae latinae} and {\em Disputationes dialecticae}. This stood in contrast to the  thirteenth century conception of Cicero's flight from public life in order to devote himself to speculation and philosophy. Now it was understood that, to Cicero, philosophy implied the possesion of technical skills and an active stance in public life: {\em Hence, as Salutati observed in his {\em De Nobilitate [ch. 23]}, the good is to be more highly valued than the merely true, virtue preferred over knowlegde, the will over the intellect, and thus rhetoric over philosophy and eloquence over wisdom.}\footnote{T.M. Conley, {\em o.c.}, p. 114} Quintilian's book made it possible to clarify those relationships, because he gave a comprehensive, coherent curriculum and subsumed all knowledge to rhetoric, which embraces everything that concerns human life. 

Only a few decades later, Manuzio's legendary Aldine Press began to publish the Greek rhetoricians: the 1508 edition of {\em Rhetores graeci} was soon followed by Demosthenes, Hermogenes,  and the other Attic rhetors. Another influential source is the Byzantine scholar Trebizond (Trapezuntius, Crete $^\circ$1395), whose most important contribution to the history of rhetoric is his {\em Rhetoricorum libri quinque}. Against Quintilian and Cicero, he holds that rhetoric should not necessarily be practiced by a ``good'' man, it is a pragmatic political art indifferent to morality. He based his book on Hermogenes's {\em On Ideas}. But Trebizond's system of mnemotechnic ``common places'' was so complicated, that he set the stage for later reforms by Agricola and Ramus, whose wide-ranging reformation of rhetorical arts started initially as an attempt to simplify the curriculum:

\begin{quote} {\em Agricola was (...) not just a reformer of the curriculum and a precursor of Ramus (...). He was at the center of a semantic revolution that reinaugurated the Classical view of invention as fundamentally rhetorical, breaking with the Boethian tradition, which restricted {\em common places} (as distinguished from particular ones) to dialectic, and returning to the Ciceronian view that had established invention and common places primarily in rhetoric and later extended them to dialectic. Agricola's critique of the tradition that separated rhetoric and dialectic parallels that of Valla's in the {\em disputationes dialecticae}, but Agricola goes further by redefining ``dialectic'' not as an art designed, as it was in the scholastic tradition, to secure the validity of an argument, but rather as the art of inquiry itself --- that is, a dialectic which is at its core rhetorical.\footnote{T.M. Conley, {\em o.c.}, p. 125.}}\end{quote}

\noindent Rudolph Agricola thus remained within the controversialist tradition, but now Ciceronian style --- developing arguments {\em in utramque partem} ---, and emancipated from its Peripathetic metaphysical foundation and its dialectical pretenses for ``truth''. His use of non-figurative ``common places'' (genus, species, quality...) as both a tool for analysis of an argument and an aid to memory explains why he could assimilate rhetorics to dialectics,  while claiming at the same time that the aim of dialectics is persuasion by probable discourse, {\em oratio}. His student Pierre de la Ram\'{e}e (Peter Ramus) radicalised these ideas in an almost opposite direction. For Ramus the aim of dialectic is again truth, {\em ratio}, but now in a Platonising context of access to truth not through controversies but through the correct, hierarchical application of the common places in the chain of an argument: {\em Indeed, properly conducted, dialectical inquiry (....) puts us in contact with ``all the multitude of things as they are in GodÕs own mind'' (dial. Inst., 38f).}\footnote{T.M. Conley, {\em o.c.}, p. 130.} These divisional schemes largely coincide with the mnemotechnical Ramist {\em epitome}, wherein every subject was ``arranged'' going from the general to the special through conceptual dichotomy.\footnote{F.A. Yates, {\em o.c.}, p. 232.} This goes back directly on Plato's conception of genuine dialectics as dichotomous division of concepts, exposed in the dialogues {\em Sophist, Statesman}, and {\em Philebus}.\footnote{I explained how and why Plato develops this method in a paper discussing the relation between Plato's logic and his ontology: K. Verelst, ÒOn what Ontology Is and not-IsÓ, Foundations of Science, {\bf 13}, 3, 2008.} Evidently, for Ramus rhetoric is inferior to dialectic, since it only deals with style and delivery of unscientific (not truth-related) discourse. Invention, disposition (the proper arrangment of arguments) and memory are the privilige of the inhabitants of the province of dialectics alone. Now all these seemingly very technical discussions had strong reverberations in the ongoing religious debates between moderate catholic intellectuals and the proponents of protestant Reformation, more specificially the debate between Erasmus and Luther. As a young student, Erasmus had been acquainted with Agricola, who inculcated him with an erudite sense for good style and an understanding of the ultimately undecidable nature of rational argumentation.\footnote{The extend to which authors like Erasmus, Montaigne... can be taken to be sceptical thinkers, and their relation to a strand of Renaisance philosophy called Pyrrhonism, has been studied in great detail by R.H. Popkin, {\em The History of Scepticism from Erasmus to Spinoza}, University of California Press, Berkeley \&., 1979.} Ramus exerted on the other hand considerable influence on Luther with his assertive dialectics of certainty and his disdain for rhetorical controversy, mainly through the protestant theologist Melanchton. 

As if things were not yet complicated enough, in 1555 a  bombshell was dropped in the middle of the debate by the publication in Italy of a by then virtually unknown Greek manuscript with a sensational content, the {\em peri hypsous} or {\em De sublimitate}.\footnote{Originally ascribed --- erroneously --- to the Hellenist writer Longinus. Its authorship remains unknown, but it became custom to refer to him is as ``pseudo-Longinus''. The history of its rediscovery is exposed in detail in M. Fumaroli, ``Rh\'{e}torique d'\'{e}cole et rh\'{e}torique adulte: la r\'{e}ception europ\'{e}enne du Trait\'{e} du sublime au XVIe et au XVIIe si\`{e}cle'', {\em H\'{e}ros}, pp. 388-393.} In it, its author, known as pseudo-Longinus, discusses the ultimate source on which all invention needs to feed in order to get access to truly new, original  ideas,  the {\em sublime (to hypsos)}. One can practice to gain access to this source, in spite of its divine nature, by developing the capacity to master the most vehement emotions and express them in an appropriately grand style. It is as if the orator develops a kind of potential participation in the divine breath into an actual capacity. pseudo-Longinus calls this exercise explicitly a {\em method} [{\em De Subl.}, II.3(3-4)]\footnote{Right from the start of his treatise, ps.-Longinus stresses the importance of the acquisition of the correct techniques to handle the sublime. He calls those who venture unprepared into this realm  ``a-methodon'', unqualified.  I used the Bud\'{e} critical edition:  H. Leb\`{e}gue (ed.), {\em Du sublime}, Les Belles Lettres, Paris, 1952.}, and he who masters this art has a literally unlimited source of inspiration at his disposal. He will be able to express everything in the appriopriate way immediately, like a prophet or a seer. ps.-Longinus compares its effect with that of a bold of lightning [I.4(9-10)]. One of the main instruments to achieve this goal is by the intensive training of rhetorical memory. 

So, during Renaissance, these ancient memory techniques made a forceful re-appearance on the public scene, but now in the potentially explosive mixture of neo-Platonic natural magic, Hermetism, and their concomittant urge to expand the powers of the mind. Interest in these mental powers had been resurrected already before, by the translation and publication by Marsilio Ficino of a Greek manuscript containing the main part of the {\em Corpus Hermeticum}, ordered by Cosimo deÕ Medici and published in 1464. Thus the essential trait of Renaissance mnemotechnical art was that it became an undisentangible alloy with religious enthousiasm, which caused a huge and increasingly hostile respons. 

 \begin{quote} {\em
Cette ``conqu\`{e}te mystique'', qui avait eu plus t\^{o}t son pendant en Italie et en Espagne, fut contemporaine, dans toute l'Europe, d'une v\'{e}ritable obsession de d\'{e}monologie et de sorcellerie, et dans les pays protestants, d'une prolif\'{e}ration des sectes d'enthousiastes.}\footnote{M. Fumaroli, ``Le cr\'{e}puscule de l'enthousias\-me'', in: {\em H\'{e}ros et Orateurs,} Droz, Gen\`{e}ve, 1996, p. 349.} \end{quote}

\noindent  A curious and interesting example of this alloy of mystical inspiration and rhetorical mnemotechnical method is Ioannis Mauburnus's (Jan Mombaer) {\em Rosetum. Excercitiorum Spiritualium et Sacrarum}.\footnote{The first edition appeared in 1491, in Basle. The work is discussed   in E. Benz, {\em Meditation, Musik und Tanz. \"{U}ber de ``Handpsalter'', eine sp\"{a}tmittelalterliche Meditationsform uas dem Rosetum des Mauburnus}, Akademie der Wissenschaften und der Literatur Mainz, Steiner Verlag, Wiesbaden, 1976. He also gives the drawing reproduced below. A detailed analysis of its content in K. Pantsers, {\em De kardinale deugden in de Lage Landen}, Uitgeverij Verloren, 2007, pp. 187-197.} Mombaer, born in Brussels, was a friend of Tho\-mas a Kempis.  He belonged to a reputed Mediaeval school of Rhineland-Flemish mysticism (Ruusbroec, Eckhardt), which was the source of inspiration for  the Dutch theologist Geert Groote when developing his {\em devotio moderna}.\footnote{B. Todoroff, {\em Laat heb ik je liefgehad. Christelijke mystiek van Jezus tot nu}, Davidsfonds, Leuven, 2002, p. 227 sq.} This seems innocent enough. But Agricola  as a young boy was educated by the Brethren of the Common Life, a religious order that belonged to this influential strand of mystical Christianity. He moreover studied in the Augustinian monastery at Erfurt, where also Luther had learned to meditate by using the {\em Rosetum}. Indeed, even Luther has known the book,  for he refers to it in his 1513 {\em Psalmenvorlesungen}.\footnote{M. Nicol, {\em Meditation bei Luther}, Vandenhoeck \& Ruprecht, G\"{o}ttingen, 1990, p. 40 sq.} Finally,  Erasmus was personally acquainted with Mauburnus, witness the correspondence between the two.\footnote{E.F. Rice, ``Erasmus and the Religious Tradition 1495-1499'', {\em Journal of the History of Ideas}, {\bf 1}, 4, 1950.} Mombaer gives an overview of the Christian mystical tradition from Antiquity up to Pico de la Mirandola, and describes a technique which uses the human hand as the ``place'' to organise memory as an aid to spiritual meditation (but behold the columns which line up the drawing at the sides of the page).\footnote{Evidently I didn't know Carruthers's 2008 article when I presented the first draft of this paper at the 2006 Cambridge conference, but I'm happy to note that my hesitating attempts at establishing this connection are confirmed by the work of such a distinguished specialist in the field. Cfr. the examples she discusses on pp. 12-13 of her {\em Mechanisms}-paper.} Mombaer's book was one of the first ever to be printed in Europe and it had an enormous succes. Ironically enough, the theological reaction championed by precisely Augustinian catholicism (like Jansenism) and protestantism against this kind of spiritual communion with the Deity would become more and more outspoken in the course of the XVI-th century. We shall come back to this below.

\begin{center}
\includegraphics[width=.7\textwidth]{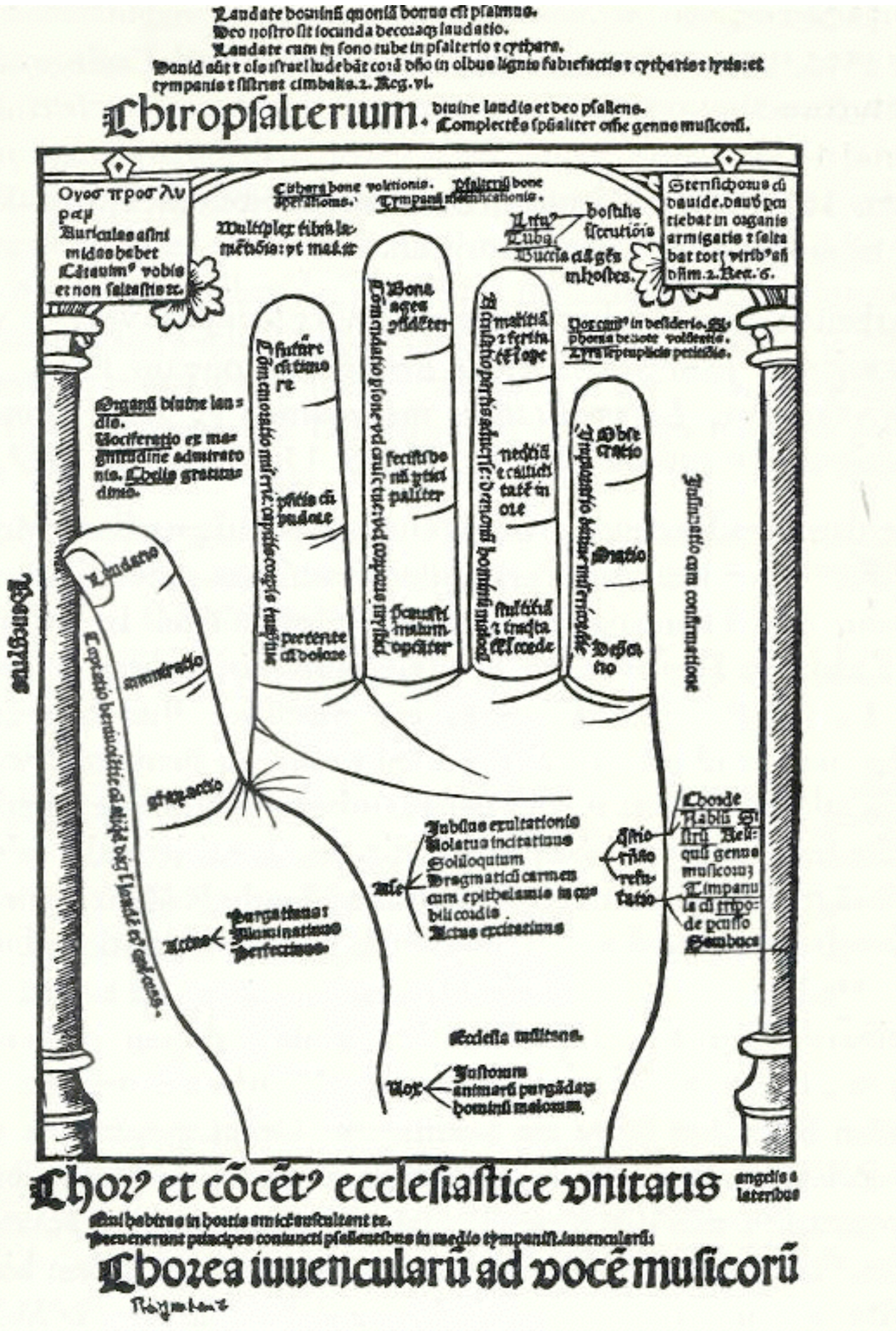}\\
\vspace{-2cm} 
{Mombaer's Handpsalter, from his {\em Rosetum}, Basle, 1491.}
\end{center}
%\vspace{-1.5cm}

\noindent But these softly Platonizing tendencies stayed carefully within the boundaries of Christianity. There were those who called for a much more radical philosophical reform.  The work exemplary for this radical ``Hermetic'' Renaissance philosophy is that of Giordano Bruno, who presents us with the baffling mixture of a superbly mathematically skilled {\em magus} who holds very modern conceptions on infinity, while practicing at the same time the most impenetrable astrological arts, and devoting himself to practising the intensification of imagination by the use of rhetoric memory as a magical tool.\footnote{See the brilliant study by F. A. Yates, {\em Giordano Bruno and the Hermetic Tradition}, The University of Chicago Press, Chicago and London, 1991 [1964]. For a point of view critical to her approach, H. Gatti, {\em Giordano Bruno and Renaissance Science}, Cornell University Press, Ithaca and London, 1999. The Lullian wheel reprinted here is discussed on p. 182 of Yates's {\em Memory}-book.} 

 One of Bruno's aims was to `formalise' the art of memory by using Ramon Lull's ``combinatorial wheels'' with the traditional technique of enhanced inner representation.\footnote{For a detailed study of the first stage of Bruno's Lullism, see M. Mertens, ``A Perspective on Bruno's  {\em De Compendiosa Architectura et Complemento Artis Lulli}'', {\em Bruniana \& Campanalliana}, iv, 2, 2009.} Lull, a contemporary to Thomas Aquinas, developed an astonishing ``Art'', the aim of which was to organise all knowledge on the basis of a meditation on the {\bf names of God} (the {\em Dignitates Dei}). There are nine fundamental names of God organised in triads (goodness, greatness, eternity,...), which form a kind of ladder of being along which the `artista' ascends and descends. At the same time they ground a memory system organised by means of combinatory wheels --- labeled from $B$ to $K$ --- laying out the interconnections that everything in creation entertains with the attributes of the Creator: 

\begin{center}
\includegraphics[width=.6\textwidth]{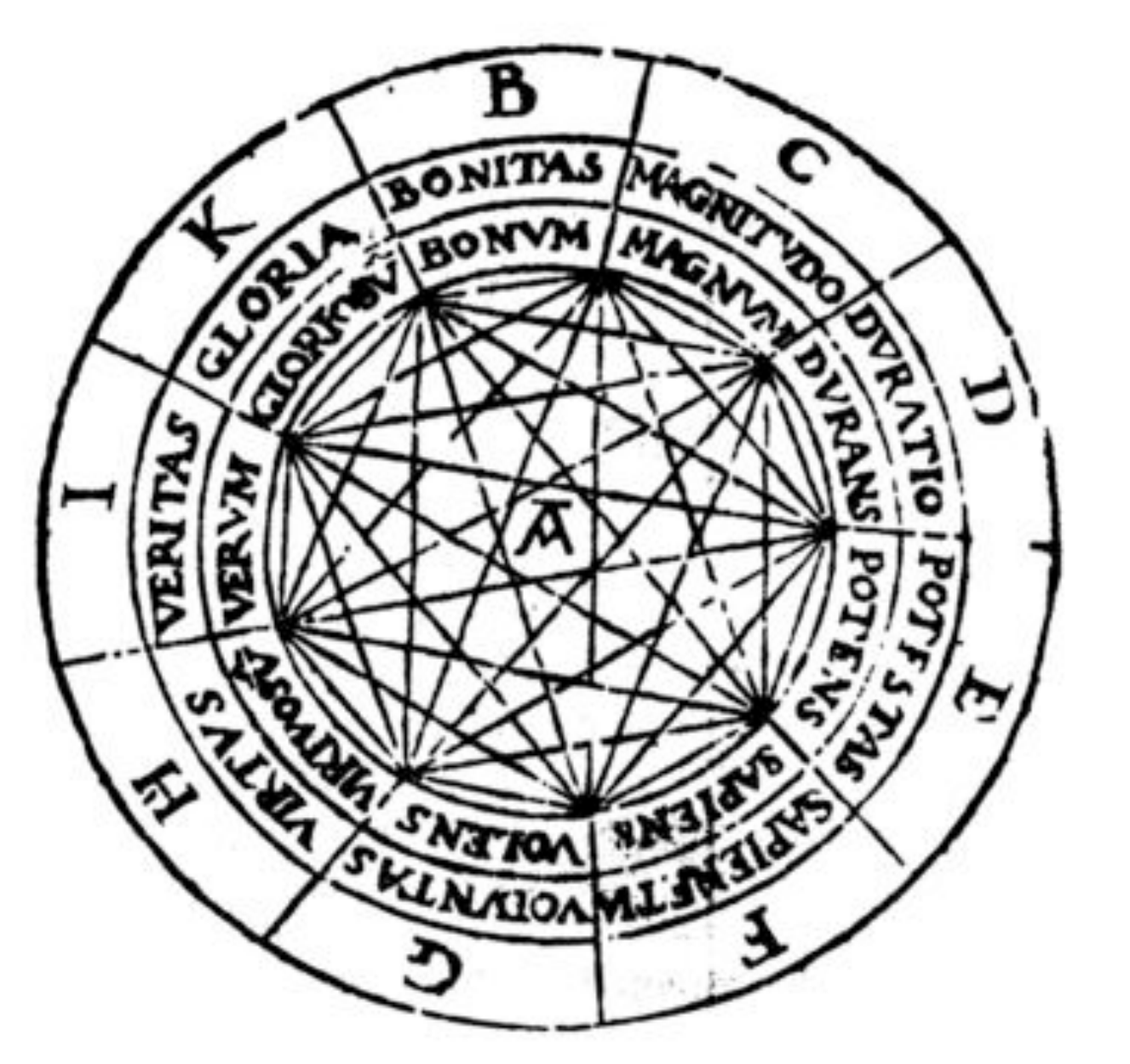} \\
%\vspace{-3cm} 
{A Lullian wheel, from his {\em Ars brevis}, in:  {\em Opera}, Strasburg, 1617.}
\end{center}

\noindent Yates gives the following explanation: {\em The geometry of the elemental structure of the world of nature combines with the divine structure of its issue out of the Divine names to form the Universal Art which can be used on all subjects because the mind works through it with a logic which is patterned on the universe.}\footnote{F.A. Yates, {\em o.c.}, p. 181.} It is a system which is universally valid and that allows the initiated to `see' at once the multiple realisations of that universal truth: He claimed that a universally valid method for attaining knowledge on the whole of the `cosmic building' was within his reach. All this philosophical brilliance aimed at restoring the direct connection between man and his Creator, and was really part of a pan-European neo-Platonic movement which had an enormous influence and which lies at the heart of the --- literal --- sublimity of the cultural heritage of European High Renaissance.\footnote{ I.P. Couliano, {\em Eros et Magie \`{a} la Renaissance}, Flammarion, Paris, 1984.}

\newpage
\subsection*{\sc The Reformation of Memory}
 
\noindent A critical response --- as well as a growing number of accusations of heresy --- was making its way. A first line of attack was opened already before the publication of {\em  On the Sublime}, by an influential Aristotelian at Leyde, J.-C. Scaliger, who wanted to counter Cardanus's {\em De Subtilitate}, in which the existence of {\em splendor} a kind of direct mystical knowledge was defended. Scaliger claims that the relation between objects and our knowledge of them is the same as that between knowledge and its expression by rational means, and that this is all there is to it. A second, more devastating line of attack was launched by Isaac Casaubon, who destroyed the historical basis of Ficino's neo-Platonic natural philosophy, by demonstrating in 1614 on philological grounds that the {\em Corpus Hermeticum} dated from the Hellenistic period, and not, as Ficino believed, from times immemorial. The in many respect decisive move was made by his son Meric Casaubon, a professor at Oxford, in 1655, when he published his {\em Treatise concerning Enthousiasm}. He dealt a final blow to the credibility of whatever kind of {\em furor} by relocating it into the realm of the pathological.\footnote{He plays therefore a key r\^{o}le in the shift from bodily to mental normality, as it has been described by Foucault in his {\em l'Histoire de la folie \`{a} l'\^{a}ge classique}.} This line of thought was taken up again by his Cambridge colleague Henry More, a Platonist of the rationalist kind, in his book {\em A brief discourse of the nature, causes, kinds and cure of Enthusiasm}.\footnote{It is well known that Newton as a young student was an avid reader of his work. J. Henry, ``Henry More and Newton's GravityÓ, {\em History of Science}, {\bf 31}, 1993. It is relevant to our concerns to point out that More discusses how God's {\em tituli} precisely serve as an indication of the abyssal gap between Him and his creation. This discussion clearly influenced Newton when writing his General scholium to the {\em Principia}. For an analysis of the textual relationshiop, see: R. De Smet and  K. Verelst, "NewtonÕs Scholium Generale. The Platonic and Stoic legacy: Philo, Justus Lipsius and the Cambridge Platonists", {\em History of 
Science}, xxxix, 2001, pp. 10-11.} 

In about the same period, another line of attack was opened by Galilei, who ridiculed the scholars that used the works of Aristotle as a kind of store-house from which they could extract pieces of textual evidence and combine them is a way as to provide an answer to whatever problem one put on the table. Thus, F. Hallyn remarks, the texts of Aristotle constitute in themselves a topical memory system: {\em une veritable {\em topique} (...) qui contient en soi tous les lieux de l'{\em invention} scientifique}.\footnote{F. Hallyn, Dialectique et rh\'{e}torique devant la ``nouvelle science'' du VXIIe sci\`{e}cle'', in: M. Fumaroli, {\em Histoire de la rh\'{e}torique dans l'Europe moderne}, PUF, Paris, 1999, p. 606.} Scientists following this `method' to arrive at new ideas are, according to Galilei in a for more than one reason remarkable comparison, like those painters who put very strict but completely arbitrary constraints on themselves while composing a portrait of a person by means of vegetables only. The crucial step in Galilei's (evidently itself rhetorical) argument occurs where he says that it is as if the ``book of the world'' was written only for Aristotle to be seen with his own eyes, and for us to rely on him for all posterity.\footnote{In his third letter on sun spots, see ``Istoria e demonstrazioni itorno alla macchie solari, terza littera'' (1613), in:  {\em Opere de Galileo Galilei}, v, p. 190.} Galilei proposes to replace Aristotle by the Book of Nature as the {\em Topica} for true invention about the world, lying open for eveyone to see and willing to study it.\footnote{F. Hallyn, {\em o.c.}, p. 607.} This appropriation of the metaphor of the ``Book of Nature''  --- which returns in the {\em dialogue on the Two World Systems} --- is al the more remarkable, because it had been used in the past over and over again to make the opposite point, namely, that scripture and nature refer to the same thing, because they stand in an allegorical relation that ultimately points back to the origin of everything, God. Augustine had started this tradition [{\em Confessions}, XIII.xv], followed by Hugh of St. Victor [{\em De tribus diebus} 4]: {\em For the whole sensible world is like a kind of book written by the finger of God (...) each particular creature is somewhat like a figure, not invented by human decision, but instituted by the divine will to manifest the invisible things of GodÕs wisdom.}\footnote{The quote is in P. Harrison's public lecture, delivered at Cambridge, May 24th, 2005. The text can be downloaded at http://www.st-edmunds.cam.ac.uk/CIS/Harrison/.} This approach climaxed in the work of the Catalan theologist Raymond Sebond, who wrote around 1435 a {\em Natural Theology or Book of Creatures}. It was translated and commented upon by Michel de Montaigne in an {\em Apology for Raymond Sebond}, and published as chapter 12 in the second part of his {\em Essays} as late as 1580! So, exactly as is the case with Descartes's occupation of the territory of `method' for the ``new science'', Galilei's conquest of the `Book of Nature' reveals as much as destroys the tradition on the borderline of which they stand, and which will be crossed irreversably only by their successors.

But before we can make truly sense of these critcisms, we need to better understand how the traditional rhetorical memory works. And in order to understand rhetorical memory, we need to take a closer look at ancient psychology, which diverges in important aspects from the modern understanding of the functioning of the human mind. The ancients hold that thoughts are made of images, and these images constitute themselves the basic material for {\em cogitatio}, active thinking. But these images are by no means direct representations of bodily sense impressions; they are not mere mental copies of the things perceived.  Sense impressions are collected in the {\em sensus communis} or {\em phantasia}; afterwards they are processed by the ``image-forming ability, {\em imaginatio} or {\em vis formalis}''.\footnote{M. Carruthers, {\em Mechanisms}, p. 4.} So an experience becomes an image through an active power ({\em vis}) of the mind.  The result of the combined activity of both organs (located in the brain) are mental ``images having formal properties that are perceptible and useful to human thought'', i.e., they are accessible to development in argumentation. Moreover, impressions are always subject to an act of estimation, a kind of instinctive judgment through yet another mental power, the {\em vis aestimativa}: 

\begin{quote} {\em Conceptions ({\em imagines}) are thus constructed by the mind from all the materials of sensation, and they have two characteristics; likeness or {\em similitudo}, and also a ``feeling'' that marks them emotionally. Thus, in this psychology there is no such thing as a completely neutral or ``objective'' experience, since all the images or {\em phantasiai} in which we comprehend our experience are already colored with some feeling before we can ``know'' them.\footnote{M. Carruthers, {\em Mechanisms}, p. 5.}}\end{quote}

\noindent Thus, the acquisition of the building blocks for the production of ideas in the mind is far from a passive process! Mental images can be stored in and retrieved again from {\em memoria} by {\em cogitatio}; the active process of recollection taking place through a valve in the brain which connects the two functional realms: the {\em vermis}. Storage should be such that images can be retraced without too much difficulty; they have to be marked and arranged in the mental places ({\em topica} or {\em loci}) that make up our memory.  Ancient memory is spatial, not logical by nature. Now, how does this rhetorical memory work? Here, again, the ancient sources at our disposal are Aristotle's {\em Rhetorica}; Cicero's {\em De Inventione} and {\em De Oratore}; the anonymous tract {\em Ad Herennium}, and pseudo-Longinus's {\em De Sublimitate}. These sources are mainly Roman or Hellenistic, but  testimonia going back to the pre-Socratic period exist.\footnote{Pliny the Elder tells us that Metrodous of Scepsis perfected the memory that was invented by Simonides of Ceos [{\em Natural History}, 7.24]. This Simonides was a fifth century lyric poet.} The baseline is that, in order to remember things afterwards (according to Aristtole: memory is a thing of the past), one has to arrange visually strong mental representations in a internalised spatial set-up. The traditional way was to internalise the complete structure of architectural artefacts like houses, palaces, squares etc., and to fill these places ({\bf \em loci}) according to the logic of their groundplan with vivid images ({\bf \em imagines agentes}) that represent the things to be remembered. As we have seen, such images are both rationally and emotionally structured; the capacity to be an agent (to `touch' the soul) is linked to their emotional content, not just to the idea of bodily movement. The author of {\em Ad Herennium} stresses the fact that places should be bright, not shadowy, and not be taken from crowed areas, for passing-by people might confuse the searching mind! [III, xiv, 31]!  Images should not be vague or banal, but exceptional in their beauty or singular in their uglines, wearing crowns or purple cloaks, or being disfigured or moving [III. xxii, 37].\footnote{I used the edition in the Loeb Classical Library: E.H. Warmington (ed.), {\em Cicero. Ad C. Herennium Libri IV}, Harvard University Press, Cambridge, Mass., 1968 [1954].}  Since a classical building always follows fixed architectural laws, this method is truly a {\em method}: the logic of the building determines the logic of recollection. Beware, however, that it is {\em not} the groundplan but the building itself that should be erected in the mind: its vivid quality is an essential ingredient to the succesful application of the method:

\begin{quote} {\em
He (Simonides) inferred that persons desiring to train this faculty (of memory) must select places and  form mental images of the things they wish to remember and store those images in the places, so that the order of the places will preserve the order of the things, and the images of the things will denote the things themselves, and we shall employ the places and images respectively as a wax writing-tablet and the letters written on it.}\footnote{F.A. Yates, {\em o.c.}, p. 2 (the quote stems from Cicero's {\em De oratore}).} \end{quote}

\noindent Penetrating internal visions make words and ideas come back ``automatically'' when one re-enters the memory-storehouse. After using them, one could empty the places of the inner architecture and re-use them afterwards. Or one could keep them and build another mental palace on top of it. Astonishing testimonies of the mental capacities thus developed abound. Augustine speaks literally about ``infinite memory''. This explains why {\em inventio} was also largely dependent on the capacity to collect and arrange common or particular places in the inner architecture of one's mind. This use of memory extends throughout the Middle Ages and well into early Modernity: {\em One should think of a properly inventoried medieval scholar's memory store as distributed across and linked throughout an organized network, a network that can be entered at any of {\bf \em  an infinite number} of places by means of a mark ({\em nota}) that is also a key ({\em clavis}), providing access from any point to the rich treasury of memory.}\footnote{M. Carruthers, {\em o.c.}, p. 7 (my bold).} 

\vspace{-0.5cm}
\begin{center}
\includegraphics[width=.7\textwidth]{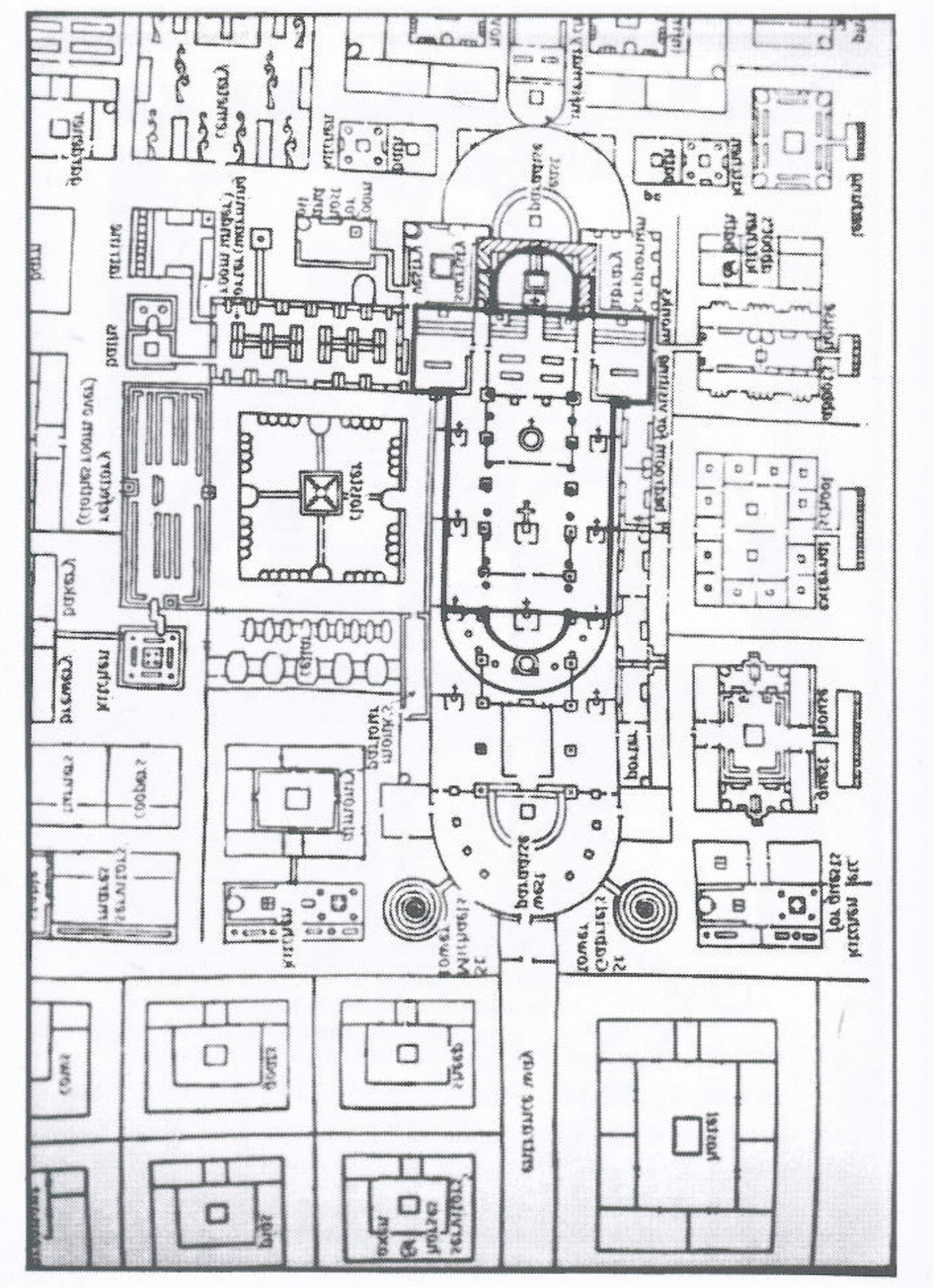}\\
\vspace{-2.5cm} 
{The monastery of St. Gall memory palace}
\end{center}

\noindent In the Mediaeval period it becomes customary to use plans instead of the buildings themselves as recollection and meditation devices. The maps may well be completely imaginary, idealised arrangements of some utopian arcchitecture like the example reprinted above.\footnote{Carruthers, {\em Mechanisms}, p. 15.} This idealistic feature brings out the methodological backbone of the practice, for it stresses the idea of a preferred itinerary to be followed,--- the moral path through monastic life mirrored in the quiet journey through the monastery's houses and gardens.\footnote{The pathway-aspect of the ``meta hodos'' has been discussed by Derrida with respect to Descartes's {\em Discours}... See M.H. Keefer, {\em o.c.}, pp. 30-36.} 

This more `abstract' use of the place-memory maybe foreshadows to some extent the later Ramist reformation. Ramus himself, however, was inspired by Quintilian, who had criticised the Greek defenders of rhetorical memory. Quintilian proposed a ``dialectical method'' to organise memory instead, based on ``dividing and composing'' ({\em Inst. Orat.}, XI, ii, 36).  As already pointed out, these divisional schemes largely coincide with the mnemotechnical Ramist {\em epitome}, wherein every subject was ``arranged'' going from the general to the special through conceptual dichotomy.\footnote{F.A. Yates, {\em o.c.}, p. 232.}  In fact, with Ramus, logic and memory coincide, but it is with an abstract, imageless memory that calls upon the use of intellectual powers only, not on those of imagination: {\em The `natural' stimulus for memory is now not the emotionally exciting memory image; it is the abstract order of dialectical analysis, which is yet, for Ramus, natural, since dialectical order is natural to the mind.}\footnote{F.A. Yates, {\em o.c.}, p. 234.} This ideal of simplicity and naturalness, as supposed to the sophistication of the artificiality of traditional memory, is a common theme underscoring the different reformatory parties, both protestant and catholic.

\begin{quote} {\em The Renaissance up to the end of the sixteenth century is characterised by an emerging conflict between the humanism of Ciceronian {\em controversia} and the newer expressions of a platonizing tendency to restore dialectic as the fundamental intellectual discipline. This conflict comes out clearly in the debates between Erasmus and Luther, which pit a view of rhetoric that is basically bilateral and symmetrical and that aims at the {\em satis probabile}, against a view like Plato's, where the Truth is to be communicated effectively in a unilateral, assymetrical setting. (...) Eventually, the humanistic controversialist rhetoric gave way, first, to the platonizing Ramists and, afterward (...) to the empiricistsÕ and CartesiansÕ (...). \footnote{T.M. Conley, {\em o.c.}, pp. 143-144.}}\end{quote}
 
\noindent Especially access through memory to the (Lullian) {\em dignitates Dei} or names or properties of God under the form of symbolic {\em images}, created a forceful reaction from calvinist and puritan sides. They called upon Ramus's imageless dialectical theory of memory as an alternative for the mnemotecnical ``idolatry''.  The puritan theologian W. Perkins even called for an absolute prohibition of the use of images as an aid to the mind! But also catholics who stood under the influence of Augustinianism, like the Jansenists, were very critical. Their common concern is the desire to purify Christianity from all remnants of paganism: {\em C'est sous l'influence d'un augustinisme plus soucieux de trancher entre paganisme et christianisme, entre nature et gr\^{a}ce, que (...) la notion d'inspiration enthousiaste (...) se verra menac\'{e}e m\`{e}me dans le domaine rh\'{e}torique et po\'{e}tique.}\footnote{M. Fumaroli, ``Cr\'{e}puscule'', in {\em H\'{e}ros}, p. 351.} 

In the introduction to this paper I mentioned Arnauld and Nicole, authors of the {\em Port Royal Logic}. The orginal title of the book is: {\em La {\bf \em Logique}, ou l'{\bf \em art de penser}}. It evidently fits in the reform-movement we outlined before, together with other attempts to revolutionise the ``art of thought'', like Johannes Claubergs {\em Logica Vetus et Nova} (1654), which might have been one of its sources of inspiration.\footnote{See the introduction to J. Clauberg, {\em Logique Ancienne et Nouvelle}, J. Lagr\'{e}e and G. Coqui (eds.), Vrin, Paris, 2007.} Arnauld \& Nicole's principal innovation is the distinction between  a thing and its representation, like a word and the  thought that it expresses: {\em (...) les mots sont des sons distinctes \& articul\'{e}s, dont les hommes ont fait des signes pour marquer ce qui se passe dans leur esprit} [II.1].\footnote{{\em Logique}, p. 143.} This seemingly innocent statement marks the fundamental breach with the past: even though disposition, style, and delivery are method-governed disciplines, the characteristic proper to inspired discourse is that the word and the thing it represents disappears: {\em Le sublime (...) abolit la b\'{e}ance entre pr\'{e}sence at repr\'{e}sentation}.\footnote{M. Fumaroli, ``R\'{e}ception'', in: {\em H\'{e}ros}, p. 383.}As Foucault points out appropriately\footnote{M. Foucault, {\em Les mots et les choses. Une archŽologie des sciences humaines}, Gallimard, Paris, 1966, ch. II.}, this separation of words and things marks out the transition from the XVI-th to the XVII-th century, and will be part and parcel of any modern liguistic theory: {\em Ainsi le signe enferme deux id\'{e}es:  l'une de la chose qui repr\'{e}sente; l'autre de la chose repr\'{e}sent\'{e}e} [I.4].\footnote{{\em Logique}, p. 80.}  But Arnauld and Nicole operate their distinction in a much larger sense, so as to separe cause and effect in natural phaenomena, or even in theology, like the symptomes of an illness, ot the stigmata of the suffering Christ [I.4]. A judgment concerning something is then a proposition consisting of this subject and a praedicate. Starting from this analysis, they defend a new method for invention based on logic to deduce new ideas. The conception of logic they defend is anti-Aristotelian in the sense that it is based on mathematical deductions as opposed to syllogistic ones. They evidently have geometry in mind, which remains after all, as far as proofs are concerned, ``an exercise in logic, classically''.\footnote{D. Finkelstein, ``Matter, Space and Logic'', in: {\em The Logico-algebraic Approach to Quantum Mechanics II. Boston Studies in the Philosophy of Science, Proceedings of the Boston colloquium for the philosophy of science V}, C.A. Hooker (ed.), 1966, pp. 199-215.} Why then geometry, instead of pure logic? Because of a problem that had been pointed out, among others, already by Descartes, and which has been summed up nicely by J. Henry: {\em The trouble with 'reason', as was clear from the fact that both Calvinists and Papists could claim it for their side, was that it could be made to subserve virtually any cause}.\footnote{J. Henry,``The Scientific Revolution in England'', in: {\em The Scientific Revolution in National Context}, R. Porter and M. Teich (eds.), Cambridge University Press, Cambridge, 1992, p. 193.} Geometry, by virtue of the fact that it has a proper content related to facts of nature, seemed to be a more reliable guided to certainty than the subtleties of controversial disputation. Arnauld and Nicole then, predictably, proceed in Part III to a critical discussion of the ``lieux'', the traditional {\em loci} or common places that were in use to structure rhetorical memory. The remarkable thing, however, is the explicit influence exerted by Arnauld's religious affiliations. Already the first page of the first chapter of the {\em Logique} contains a reference to St. Augustine. And he was a theologist adhering to the ``movement of Port Royal'', a catholic heresy also known under the name {\em Jansenism}, for its founder, Cornelius Jansen, theologist at Louvain and bishop of Ypres, who published (posthumously) a book on Augustine in 1640. From this book five heretical theses were extracted, and condamnation followed in 1643.\footnote{J.P. Chantin, {\em le jans\'{e}nisme}, Paris, Cerf, 2000.} Jansenist theology is Augustinian by inspiration and marked out by strong resemblances to several tendencies in radical protestant Enlightment, i.e., Calvinist, puritan, millenarian and prophetic. It included positions on the role of God's grace, not ot be derserved by everybody but only for a group of `elected', with the ensuing obligation to moral purity, and the remarkable but related thesis that Jesus the Saviour was not born to save everybody; strong tenets on praedestination; a belief in the prophetic renewal of the Church, etc. The French philosopher Pascal belonged to their followers as well. It was vehemently anti-papal, and for quite a while on its way to constitute the intellectual backbone of an ``\'{E}glise gallicaine'', a strand of French catholicism much more independent from ecclesiastical authority than the traditional ``Roman'' one. They were condemned several times after 1643, and by different popes (Alexander VII, Bull {\em Ad sacram} (1656); Clemens XI, Bull {\em Unigenitus Dei Filius} (1711).\footnote{L. Cognet, {\em Le Jans\'{e}nisme}, PUF, Paris, 1995 [1961], p. 72; p. 85 sq.} The Jansenists called for a rationalisation and purification of both language and of the {\em ars inveniendi}, not in the sense of Ramus, but in the sense of Descartes. Arnauld and Nicole are also the authors of a {\em Grammaire g\'{e}n\'{e}rale et raisonn\'{e}e contenant les fondemens de l'art de parler, expliqu\'{e}s d'une mani\`{e}re claire et naturelle}\footnote{A. Arnauld et C. Lancelot, {\em Grammaire g\'{e}n\'{e}rale et raisonn\'{e}e}, \'{E}ditions Allia, Paris, 1997 [1660].}, which essentially claims that  the basic structures of language are simple, innate and universal.\footnote{The work inspired Noam Chomsky in the development of his ``generative grammar'', which he himself qualifies as a part of a wider approach to language, ``Cartesian linguistics''. see N. Chomsky, {\em Cartesian Linguistics. A Chapter in the history of Rationalistic Thought}, New York and London, 1966.}  The separation of mentalisties can be gauged from the fact that C. Favre de Vaugelas, co-founder of the {\em Acad\'{e}mie} and one of the main editors of its {\em Dictionnaire}, published only fifteen years before a  book {\em Remarques sur la langue franaise, utiles ˆ ceux qui veulent bien parler et bien \'{e}crire}, of which the central tenet was that the best way to master the French language was to copy the ways of speech --- ``le bon usage'' --- of the members of the royal court.\\

For Meric Casaubon, the poetically transgressing ``inspired'' mind is nothing more than a pahtological mind, vulnerable to the ruses of paganism. True judgment requires nothing but ``bare speech''  to expose ``the true nature of things''. Words are nothing but representations, the mental face of sensible things. Thus a conceptual line links the in origin Aristotelian conception of words as {\em symbola rerum} defended by both Scaliger and Casaubon, directly to Locke's {\em Essay on Human Understanding} (1690). Indeed, nobody will be surprised to find that it contains a critique of enthousiasm which insists on the lack of proof for the pretenses to certain knowledge that any utterance grounded in it could claim. On the other hand, Henry More succesfully ridiculises the melancholic temperaments unable to discriminate between a bodily affliction and the workings of God's grace.\footnote{M. Fumaroli, ``Cr\'{e}puscule'', in {\em H\'{e}ros}, pp; 372-372.} While the English Aristotelians and later empirists thus succeeded in effectively abolishing any credible remnants of the ``inspired'' forms of knowledge to the extend that all activity involving possible errants into {\em l'irrationnel, l'inspiration, la dict\'{e}e ``d\'{e}monique''} --- even philosophy ---- became suspect, the alliance between ``simplicit\'{e}'' and ``le sublime'' was to shape what Ren\'{e} Bray was to call ``la doctrine classique'', that gracious philosophical ``art de converser'' in urban life in France's XVII-XVIII-th centuries. Even the detractors of rhetorics had to put their cause eloquently! During that period, the debate on the nature of the sublime will resurface, but now within the safe confines of the domain of the arts, in the modern sense of that word. The notorious {\em Querelle des Anciens et Modernes}, initiated by Montaigne, was to oppose the great French dramatists Racine and Corneille, the writers Perrault and Balzac.\footnote{M. Fumaroli, {\em Querelle}, pp. 7-10.} But even Boileau, author of a {\em De l'art po\'{e}tique}, and the first translator of the {\em De Sublimitate} into French, had to surrender at least partially to Scaliger's criticisms.\footnote{M. Fumaroli, ``Cr\'{e}puscule'', in {\em H\'{e}ros}, p. 360.} This evacuation of philosophy into the realm of literature went at the cost of its immediate connection to and relevance for reality, whether mundane or sublime: imagination became a sphere of life disconnected from ``science'' as it was practiced in laboratories or at the theorist's writing desk. 

\begin{quote} {\em Car les contraintes psychologiques et m\^{e}me physiques qu'exerca la reforme de l'\'{E}glise - du c\^{o}t\'{e} protestant comme du c\^{o}t\'{e} catholique - ne furent que de peu inf\'{e}rieures \`{a} celles exerc\'{e}es par la R\'{e}volution francaise \`{a} son apog\'{e}e ou - mutatis mutandis - par la r\'{e}volution sovi\'{e}tique. (...) A un moment donn\'{e}, la censure avait transform\'{e} la personnalit\'{e}: les gens avaient perdu l'habitude d'utiliser activement leur imagination et de penser par ``qualit\'{e}s'', car cela n'\'{e}tait plus permis. La perte de la facult\'{e} d'imagination active entra\^{\i}na forc\'{e}ment avec elle l'observation rigoureuse du monde mat\'{e}riel et celle-ci se traduisit par une attitude de respect pour toute donn\'{e}e quantitative et de soup\c{c}on envers toute assertion d'ordre ``qualitatif''.}\footnote{I.P. Couliano, {\em Eros et Magie \`{a} la Renaissance}, Flammarion, Paris, 1984, p. 240.}\end{quote}

\noindent When the gates of the soul finally closed on the sublime as a road to ultimate reality, the way was free for a strictly rational account of the functioning of the human mind in its relation to the world. Both the empirists and the Cartesio-Leibnizian rationalists are holding the same line with respect to this, be it armed with very different argumentive weaponry. Style becomes a facet of the first person in the small sense, --- the birth of the subject, so to speak, and the birth of the idea that debates are a matter of mere ``personal opinion''. The question of knowledge belongs from now on to an entirely different realm. The empirists evacuate the question of the sources of knowledge to the outside world, and reduce it to the relation between that world and the sensual impressions we receive from it. This approach will be epitomised in the {\em Regulae Philosophandi} by which the mature Newton starts Book III of the {\em Principia}. The clearest example of the nature of the periods' concern with method, however, is Leibniz, who explicitly links the trustworthliness of his {\em caracteristica universalis} to the exhaustive success of his {\em ars inveniendi}. And all are deeply indepted to Descartes's ideas. We by now have reached the point were the story of Graham and Kantor begins.

\subsection*{\sc When Existence is in a Name}

\noindent During the spring of 1913, the Russian navy attacked the monasteries on the peninsula of Athos in order to destroy a Christian-orthodox heresy that had been condemned shortly before by the synod of St. Petersburg. More than 1000 monks were arrested and imprisoned. What was the nature of this heresy, and why did it call for such a violent reaction? Its roots go back to a problem that had haunted theology since centuries: the question whether God can be known, and if not, how He can be worshipped.\footnote{This is in fact a variation on Meno's knowledge paradox, presented by Plato in the dialogue bearing the same name: how can one learn what one does not know, given the fact that one does not know one does not know it?} In the Western or Latin church, answers to this question had been provided along two different lines: that God is to a certain extend rationally knowable through His creation, but at the same time, as far as His mysterious nature is concerned, He remains out of reach for mortal minds, and one can trust only on faith, as revealed in Scripture and transmitted by the Church. The consequences of this ambivalence had stirred Christianity in its foundations already in the thirteenth century, for it developed into the {\em theory of double truth}\footnote{E. Gilson, ``La doctrine de la double v\'{e}rit\'{e},'' in: {\em \'{E}tudes de philosophie m\'{e}di\'{e}vale}, Universit\'{e} de Strasbourg, 1921.}, which held that if reason proved something to be true, while its contrary was a matter of faith, then both acquired validity. The followers of this school of thought were mainly Aristotelian masters at the Facutly of Arts at University the University of Paris. They are called ``Latin Averroists'', after the notorious Arab philosopher and commentator of Aristotle, Averroes (Ibn Rushd), and condemned in 1277 by bishop Etienne Tempier. In reaction to the radical Aristotelians, Thomas Aquinas developed his magistral {\em Summa Theologica}, in which the syntesis between Aristotelianism and Catholic dogma  was coined that would govern Latin orthodoxy for centuries to come. But a too sharp division between reason and faith was ad odds with the neo-Platonic currents feeding Christianity since late Antiquity, and on whose fertile soil mystical tendencies using music, poetry and meditation techniques, had flourished throughout the Middle Ages. The Mediaeval church had been rather lenient toward them, even though they were often teetering on heresy because of the theologically unacceptable idea that the exalted soul could somehow unify itself with the deity. One can clearly see the resonances with the theological and philosophical disputes on inspiration and invention we discussed above. The connection becomes explicit in the work of some of the main characters in the dramatic sequence of events: from Lull's ``crusade'' against Averroism via Ramus's defence of the Paris masters, we enter the Reformation period during which the chasm between revealed and discovered truth will deepen to the extend that it will eventually lead to Galilei's 1633 condemnation.\footnote{L. Bianchi, {\em Pour une histoire de la ``double v\'{e}rit\'{e}''}, Vrin, Paris, 2008. See the ``pr\'{e}face'' for an overview.} Nevertheless, at both sides of the gap there was agreement on the need to reject and even persecute religious ``enthousiasts'' of all stripes, but especially the rapidly spreading millenarian Protestant sects like the ``French Prophets''. They believed they could communicate with the deity directly through possessions, including visions, prophecy, and trances, induced by dancing and music. They belonged to a bigger group of radical Protestants, the Camisards, who had been violently oppressed since the revocation of the Edict of Nantes (1685). In 1706 a group of refugees lead by Elie Marion arrived in London, where they met again with oppressive authorities. Notwithstanding that oppression, they continued to exert lasting influence through their writings. Marion had written a treatise on prophecy, {\em Les Avertissements proph\'{e}tiques}\footnote{J.P. Chabrol, Elie Marion, le vagabond de Dieu 1678-1713, Aix-en-Provence, Edisud, 1999.}, which laid the basis for the later Shaker movement.\footnote{This Marion had been a close friend of Fatio de Duiliier, himself a close friend of Newton's at least until the mid 1690s.  See L. Verlet, {\em La malle de Newton}, Paris, Gallimard, 1993, p. 194.  Newton himself has written extensively on the intepretation of prophecy, and his method with respect to this relates closely to the methods he developed in his scientific works. This has been researched extensively in M. Mamiami, ``To twist the meaning: Newton's Regulae Philosophandi
revisited'', in J. Z. Buchwald and I. B. Cohen (eds.), {\em Isaac NewtonÕs Natural Philosophy}, 
Cambridge, Mass., 2001, pp. 3Ð14; id., ``Newton on prophecy and the Apocalyps'', in
I. B. Cohen and G. E. Smith (eds.), {\em The Cambridge Companion to Newton}, Cambridge, 2002, pp. 387Ð408. Accordong to Ben Chaim, Newton's understanding of causal explanations and experimental observation as rules that govern perceptual judgment, not as interpretation of given phaenomena, fit in a larger scheme: {\em His endeavour to explore and identify such rules followed his practical reflections on human beings as agents who belonged to God's dominion and were created to serve its divine ends. These reflections suggested, more specifically, that the aim of natural philosophy was the discovery of divine rules that instrumentally constrained and facilitated human conduct.} M. Ben Chaim, ``The discovery of natural goods: Newton's
vocation as an `experimental philosopher' '', {\em BJHS}, 2001, {\bf 34}, p. 397. There are indeed many good reasons to accept that the relation between early modern science and religion is far more complicated then generally assumed. For an overview of the Newton-related material, see S. Snobelen, `` `To discourse of God': Isaac Newton's heterodox theology and his natural philosophy'', in P.B. Wood, (ed.), {\em Science and dissent in England, 1688-1945},  Aldershot, Ashgate, 2004, pp. 39-65.} But in Britain their influence was stamped out by the combined efforts of the Anglican church and the philosophical empiricists, who had been taking over the intellectual lead from the rational Aristotelians of the former generation. We thus see that, even though inspired movements of all kinds had been present and even dominating throughout large parts of European cultural history, their public influence ceased to exist almost completely after, say, 1750, at both sides of the Channel, save for the literary credibility they retained in France.  This had immediate consequences for philosophy in its large sense. Descartes starts his {\em Regulae} with the remark that people tend to think wrongly that they can conclude validly on the basis of mere similarity (``similitudo'') between things that otherwise differ. We mentioned already that Foucault points out that the correction of this ``error'' marks the transition from ancient to modern with respect to linguistics. We showed how Arnauld and Nicole (inspired as well by Pascal) elaborate this distinction into a full-fledged theory of representation. From now on, ``likeness'', as well as any other mental operation that involves a leap of the imagination, has its place only in philosophy or literature, not anymore in science. Furthermore, in the first Book of his {\em Principiae philosophiae}, Descartes states clearly that, since we are finite beings, it would be absurd to hope we might ever understand the infinite, and it therefore be better to not think of it at all. He explicitly refers to mathematical problems, like whether there are infinite lines, and whether their halves are infinite as well; whether there are infinite numbers, and if so, whether they are even or odd, \&c. In order to solve such questions, you need to have an infinite mind.  They are incomprehensible to us to the same extend as the mysteries of the Incarnation and of the Trinity, which he discussed in the paragraph before. It is worthwhile to give the full quotation:

\begin{quote} {\em Ita nullis unquam fatigabimur disputationibus de infinito. Nam sane, cum simus finiti, absurdum esset nos aliquid de ipso determinari, atque sic ilud quasi finire ac comprehendere conari. Non igitur respondere curabimus iis, qui qaerunt an, si daretur linea infinita, ejus media pars esset etiam infinita; vel an numerus infinitus sit par anve inpar, \& talia.}\footnote{We use the Adam and Tannery edition of the {\em Oeuvres de Descartes}: AT VIII-1, I-26}\end{quote}

\noindent For us, finite minds, when confronted with something of which we can conceive no end, the only viable strategy is to consider it as {\em indefinite}, i.e., finite, but without a conceivable end. In contemporary terms: all arguments we develop in mathematics are finitary, even if they deal with `the infinite'. 

\begin{quote} {\em Nos autem illa omnia, in quibus sub aliqua consideratione nullum finem poterimus invenire, non quidem affirmabimus esse infinita, sed ut indefinita spectabimus.}\footnote{AT VIII, I-26}\end{quote}

\noindent He then gives infinite extension and the infinite divisibility of matter as examples of problems which we can only think of in ``indefinite'' terms. ``Infinity'' is a name we can apply to God alone: {\em tum ut nomen infiniti soli Deo reservemus} [AT VIII, I-27]. We can, however, {\em  know} that there are infinite things --- eternity, God --- without {\em understanding} them, because without them, it would not even be conceivable for a finite mind to have the idea of infinity at all. In a certain sense, our finitude points towards the infinity of God, making it accessible to our knowledge, but keeping it far out of reach from our comprehension. There are many other places in Descartes's work where he varies on this basic idea. {\em (...) our reason, our powers of comprehension are finite and limited, and since the will's decisions are determined by reason, we have no, so to speak, immediate access to infinity. We have to content ourselves with recognizing infinity without being able to explain it} [AT III 292; II 138].\footnote{See for this translation and an elaborate discussion: A. Drozdek, ``Descartes: Mathematics and Sacredness of Infinity'', {\em Laval thŽologique et philosophique}, {\bf 52},1, 1996.} This in fact comes down to the Aristotelian difference in degree of reality between `potential' and `actual' infinity. The cautious attitude with respect to actual infinity was shared by Huygens and even by Newton in his mature life.\footnote{Once again, Newton's attitude is peculiar and ridden with apparent inconsistencies. A case in point is the mathematical methodology of the {\em Principia}, which is deliberately based on the old geometrical methods, while he had his theory of fluxions already at his disposal. This was part of an ideological agenda to ptomote ``ancient wisdom'' as superior over the modern one; see J.E.Force, ``Newton, the `Ancients'' , and the `Moderns' '', in: {\em  Newton and Religion : Context, Nature, and Influence}, J. E. Force and R. H. Popkin  (eds.), Kluwer, Dordrecht, 1999, pp. 237-257. The implications this had for Newton's mathematics are brilliantly analysed in N. Guicciardini, {\em Reading the {\em Principia}, The Debate on Newton's Mathematical Methods for Natural Philosophy}, Cambridge University Press, Cambridge, 2003, [1999].}  Pascal's introduction of the principle of induction in his {\em Trait\'{e} du triangle arithm\'{e}tique} (1665), seems at first to jeopardise this line of events. But evidently, to establish his induction principle, Pascal uses arguments based on indefiniteness, not infinity, in the sense of Descartes.\footnote{Poincar\'{e} showed in his 1905 paper ``Les math\'{e}matiques et la logique", that complete induction thus construed involves either a {\em petitio principii} --- the set of all natural numbers ---, or a circularity --- the principle of complete induction itself. The paper has been published in several pieces and variants; I used the annotated reprint in G. Heinzmann, {\em Poincar\'{e}, Russell, Zermelo et Peano. Textes de la discussion (1906-1912) sur les fond\'{e}ments des math\'{e}matiques: des antinmies \`{a}  la pr\'{e}dicativit\'{e} }, Blanchard, Paris, 1986, pp. 11-54.} Moreover, Pascal had warned against epistemological pretenses: even though infinities exist, they remain infinitely far beyond the limits of comprehension of the human mind.\footnote{Cfr. the famous passage 72 on the {\em cirron} (mite) in the {\em Pens\'{e}es}. Arnauld and Nicole elaborate this ``problem of infinity'' along Pascalian lines in the first chapter of the fourth part (``De la m\'{e}thode'') of their {\em Logic}; {\em o.c.}, p. 359 sq. See also the introduction.} So there is much more direct evidence to Descartes's and Pascal's stance with regard to inifnity than the meager quotation on Descartes's reductionism Graham and Kantor offer on p. 36 of their book. We also see more clearly by now in what philosophical and religious context they should be placed. Be that as it may, according to them, this is the intellectual attitude that, almost three centuries later, hampered the French mathematical community in its attempts to deal with Cantor's set theory and its concomittant arithmetic of infinity.\footnote{ The following sketch is based on  L. Graham and J.-M. Kantor, {\em Naming Infinity. A True Story of Religious Mysticism and Mathematical Creativity}, Belknap Press, Cambridge, Massachusetts, 2009, pp. 25-40.} Cantor had proven in 1873 that the set of natural numbers (the integers) and the set of real numbers (the real line) had different kinds of infinite numbers of elements, and, even more, that there is an infinite number of infinite kinds. He then set out to check whether there are other infinities lying in between the two arithmetical infinities, and formulated his famous Continuum Hypothesis (CH), which basically states that this is not the case. Although he never succeeded in proving CH --- nobody ever did --- he did prove in 1884 the remarkable result that CH holds for all closed subsets of the real line. This idea of divisions of the real line had of course ancient antecedents in the paradoxes of Zeno, and had been problematised before by Bolzano in his {\em Paradoxien des Unendlichen}. Interestingly, one of the early defenders of Cantor's work in France was the mathematician and historian of mathematics and philosophy P. Tannery. He introduced in the modern literature the idea that set theory and the paradoxes appearing in it should be related to the work of Zeno.\footnote{If one renders Zeno's two fundamental paradoxes in what I would call their canonical form as 1) paradox of plurality: to consist of parts with and parts without magnitude, and 2) paradox of motion: to count the uncountable, than the family resemblance with the modern set theoretic paradoxes (Burali-Forti, Russell,...), becomes easely visible. I prsented my reading of Zeno's paradoxes in: K. Verelst,  ``ZenoÕs Paradoxes. A Cardinal Problem. I. On Zenonian
Plurality'', in: Paradox: Logical, Cognitive and Communicative Aspects.  Proceedings of the
First International Symposium of Cognition, Logic and Communication, Series: The Baltic
International Yearbook of Cognition, Logic and Communication, Vol. 1, University of Latvia
Press, Riga, 2006. For Tannery's ideas on the subject, see P. Tannery, ``Le concept scientifique du continu. Z\'{e}non d'El\'{e}e et Georg Cantor'', {\em Revue philosophique de la France et de lÕ\'{e}tranger}, {\bf 20}, 1885. One will remark that this is a philosophical, not a mathematical journal. It might be worth considering that Tannery is (together with Charles Adam) also the editor of the famous critical edition of Descartes's complete works!} The problem was pinned down accurately in 1905 by yet another prominent French mathematician, H. Lesbesgue (of the Lebesgue integral), when he entered the discussion on another set theoretic source of antinomies, Zermelo's axiom of Choice: {\em Peut-on s'assurer de l'existence d'un \^{e}tre math\'{e}matique sans le d\'{e}finir? D\'{e}finir veut toujours dire nommer une propri\'{e}t\'{e} du d\'{e}fini.} This is just another way of saying that one has to restrict mathematical notions to those which can be clearly defined and for which a consistent mental representation is present, requirements which sound familiar in an intellectual environment shaped by Descartes and Port Royal...

So even though French interest in Cantor's results was almost immediate, it was relegated largely and for a long time to the realm of philosophy. When three of the most noted French mathematicians decided already in 1880 that it would be worthwhile to translate Cantor's work, they gave the job to a Jesuit priest, because {\em [Cantor's] philosophical turn of mind will not be an obstacle for a translator who knows Kant}.\footnote{The quote is on p. 30 of the Graham-Kantor book.} In 1898, E. Picard, an important member of the French mathematical establishment, nearly dismissed R. Baire's doctoral dissertation on discontinuous functions. Baire had, be it in a very careful way, used concepts of set theory in the theory of functions. Evidently, such a speculative subject witnessed the twists of a philosophical mind: {\em L'auteur nous para\^{i}t avoir une tournure d'esprit favorable \`{a} l'\'{e}tude de ces questions qui sont \`{a}  la fronti\`{e}re de la math\'{e}matique et de la philosophie}. I take it that I do not have to explain anymore that this was not an innocent remark, nor meant as a compliment. Even Borel, who invented Borel sets (a way to divide the real line that has a specific algebraic structure) and founded measure theory, and who was as a young mathematician captivated by both the person and the work of Cantor, declared that he had been carried away by German romanticism, and took more distance: {\em We are serious people; this at least is not philosophy; a disagreement can only be due to a misunderstanding}.\footnote{L. Graham and J.-M. Kantor, ``Name Worshippers. Russian religious mystics and French rationalists: Mathematics 1900-1930'', {\em American Academy of Arts and Sciences Bulletin}, vol. LVIII, nr 3, 2005.} Cantor's ideas had indeed flourished in late XIX-th century romantic Germany, but even there Cantor had faced problems throughout his career, because of the staunch opposition to his work by an influential group of mathematicians led by L. Kronecker, who declared in explicitly theological terms, that ``God made the integers; all else is the work of man''.\footnote{E. T. Bell, {\em Men of Mathematics}, Simon and Schuster, New York, 1986, p. 477.} Both Cantor and Baire suffered a nervous breakdown as a consequence from the mental strain caused by working on such problems,  as well as the burden it put on the development of their respective academical careers. Cantor ended up in a mental asylum, and Baire committed suicide, thus complying to Meric Casaubon's curse in one of those ironical twists history seems so fond of to play on us.\\

How different the reception of these ideas was in Russia can be gauged from the remark of a proponent of the Russian school of mathematics, which contributed so much to their further elaboration: {\em Everything seems to be a daydream, playing with symbols, which however, yield great things.} This note was not written by an exalted artist, but by one of the greatest Russian mathematicians of the twentieth century, Nikolai N. Luzin.  Graham and Kantor explain this difference by the cultural influence dominant in Russia during the end of the XIX-th and the beginning of the XX-th century. We are not going to reproduce here in extenso their analysis --- I refer to the book and the papers already cited ---, but let us summarise their conclusions, and link them to the bigger historical context we outlined above. In the cultural realm of the Eastern (Greek and Russian) Orthodox churches, nothing comparable to the extinction of religious enthousiasm had happened until much later, not incidentially during the period of Russia's forced modernisation, started by the Czarist regime in the second half of the XIX-th century. The use of exalted symbolism had been common practice in Orthodox liturgy for ages, e.g.  in iconography, but especially the use of the controversial ``Jesus Prayer'', a ritual in which the worshipper {\em chants the names of Christ and God over and over again (...) until his whole body reaches a state of religious ecstasy in which even the beating of his heart, in addition to his breathing cycle, is supposedly in tune with the chanted words ``Christ'' and ``God''}.\footnote{L. Graham and J.-M. Kantor, {\em Isis}-paper, p. 68.} It had roots in practices going back to fourth century monastic hermits in the Near East, and caused already in the XIV-th century a theological controversy with the more rationalistic Byzantine tradition. Those practicing it were called ``Name Worshippers'' or ``Nominalists''; those opposing the practice ``anti-Name Worshippers'' or ``anti-Nominalists''. The conflict came to a climax in a politically highly sensitive period (the impending collaps of the Ottoman Empire and the recent Menchevik Revolution) when in 1907 Ilarion, a monk from a Russian monastery on Mt. Athos, published a very popular book entitled {\em In the Mountains of the Causasus}, in which he claimed that the fact that the faithful can reach a state of unity with God by chanting his name proves that the name of God is holy in itself, that the name of God {\em is} God. Interestingly enough, he stresses that the process required to learn to do the prayer effectively takes years and has to follow certain methodological steps. With his claims, he provoked a line that no monotheistic theology can afford to be trespassed, and condemnation followed swiftly, which led to the events cited at the beginning of this section.\footnote{A tasty report on the revolt at Mt. Athos and details on the larger historical background can be found in the first chapter of Graham and Kantor's book.}

Now this strand of exalted theology exerted a direct influence on the mathematical developments in Russia through a former mathematician turned Orthodox priest, P. Florenskii. Two of the founding members of the Moscow school of mathematics, D.F. Egorov and N.N. Luzin, who would become leading mathematicians in the twenties and thirtees with major contributions to descriptive set theory and the theory of discontinous functions, were among Florenskii's disciples. Florenskii was an adept of the Name Worshippers, and he defended his ideas on several occasions during meetings of a small circle of followers in Egorov's  appartement. He sought to bring together mathematics and religion through the concept of ``naming'', as it had been understood by the Nominalists: {\em to name something was to give birth to a new entity. (...) Humans could exercise Free Will and put in perspective mathematics and philosophy}.\footnote{L. Graham and J.-M. Kantor, {\em Isis}-paper, p. 70.} Set theory was the field par excellence wherein the mathematical power of naming came to full fruition: discontinuity and infinity became a matter of mental creativity instead of nervous breakdown. Luzin and Egorov were not witheld by the constraints that halted Lebesgue. They fully developed the concept of ``effective set'' or ``named set'', introduced by Lebesgue in 1904 as a tool to avoid the inconveniences of the Axiom of Choice. According to Graham and Kantor, the birth of descriptive set theory coincides with the moment in 1917 on which M. Suslin, a student of Luzin's, entered his office to show a mistake he discovered in a proof by Lebesgue concerning projections of Borel sets on the real line. The famous Polish topologist W. Sierpinski witnessed the scene.\footnote{L. Graham and J.-M. Kantor, {\em Isis}-paper, p. 72.} Sierpinski's own work would lead later on to another field infested by discontinuity and infinity, the theory of fractals. Luzin was convinced that properly naming a mathematical object was the key to many profound mathematical problems. We now better understand the motto with which this paper started. Thus, Luzin's work is the fulfilment of the project outlined by Florenskii: {\em the naming of sets was a mathematical act, just as the naming of God was a religious one}, because {\em the point where divine and human energy meet is `the symbol,' which is greater than itself}.\footnote{L. Graham and J.-M. Kantor, {\em Isis}-paper, p. 70.}\\

But this leaves us with a seeming paradox: In his book on the symbolic revolution in XVI--XVII-th century European mathematics, M. Serfati  has shown beyond a shadow of doubt that the ``de-rhetoricalisation'' of mathematics and its  increasing ``symbolisation'' go hand in hand.\footnote{M. Serfati, {\em La r\'{e}volution symbolique. La constitution de l'\'{e}criture symbolique math\'{e}matique}, P\'{e}tra, Paris, 2005.} It is nevertheless clear that the Russian use of symbols pertains to the tradition that favours the use of the imagination as a source of creativity, c.q., the rhetorical tradition. How can these two conclusions regarding the use of symbols be matched? We have no intention to go into this question in any detail here, but let us suggest the outlines of a putative answer. Apparently, in the Russian case there is another source of emancipation of the symbol at work, with a very different cultural origin. The emancipation of the symbol is not linked to its being purified from any pre-existing content referring to a ``chose'' in the sense of the early algebrists, as described brilliantly by Serfati in his book. We suggest that, for the Russians, symbols act like emblems, the existence of which is as much a condition for the existence of what they ``symbolise'' as it is the other way around. This is confirmed by the way Graham and Kantor understand the Russian way to resort to the use of symbols: {\em Symbolism is the use of a perceptible object or activity to represent to the mind the semblance of something which is not shown but realised by association with it.}\footnote{L. Graham and J.-M. Kantor, {\em Isis}-paper, p. 68.} Evidently, this is generally {\em not} the way we construe the operation of mathematical symbols. We do agree, however, that this is exactly what the Russians had in mind: {\em nommer, c'est avoir individu}. Symbols play a different r\^{o}le in the two approaches --- the rationalistic and the imaginative one. In the first, symbols are representations of something else, even a `no-thing'. In the second, symbols are signs that not simply ``re-place'' something, but they take place instead of the thing in the place of which they stand, a semiotic rather than a semantic relation, a co-incidence rather than an arbitrary union. The messenger in a certain sense {\em is} the master speaking out directly.\footnote{The mechanism of the coincidence between a person preaching on supply and his master has been studied in: A.R. Johnson, {\em The One and the Many in the Israelite Conception of God}, University of Wales Press, Cardiff, 1960, especially with respect to the workings of prophecy.}

\begin{quote}{{\em Il faut que ce soit quelque chose de c\'{e}leste et d'inspir\'{e} qui intervienne dans l'\'{e}loquence pour exciter les transports et les admirations qu'elle cherche. (...) Il est besoin de quelqu'autre que de l'art, afin que la sp\'{e}culation se rende sensible, et qu'elle tienne ce qu'elle a promis. {\bf \em Afin que les r\'{e}gles deviennent examples, afin que la connaissance soit action et que les paroles soient des choses.}}\footnote{J.-L. Guez de Balzac, {\em Paraphrase ou de la grande \'{e}loquence}, cited by M. Fumaroli, ``R\'{e}ception'', in: {\em H\'{e}ros}, p. 395.}}\end{quote}

\subsection*{\sc Conclusion}

The comparative study by Graham and Kantor seems to indicate that a `leap of imagination' --- in the traditional sense --- is needed in order to be able to reach out for certain kinds of formal ideas. We believe the historical research presented in this paper confirms their findings at least to a certain extend. The r\^{o}le of imagination, even when not recognised as such for want of relevant cultural reference points, sheds another light on the tenacity of the infamous ``Platonism'' which is rumoured to be so popular among practicing mahematicians, to the despair of philosophers of mathematics and logicians alike.  It might well be that such type of conviction, no matter how flawed in itself, is a condition of possibility for doing certain types of mathematical inquiry at all.  The dual also holds, but it remains a dual, not an inverse relationship.

Evidently, this paper is just a first and sketchy attempt to deal with this sort of question. In this specific case in order to be complete, a detailed study of the philosophical and rhetorical sources of orthodox Christianity is needed. Such an inquiry required knowledge of Byzantine Greek and Slavo-Baltic languages, as well as a fine knowledge of Orthodox theology, something to which I can not even aspire. But I hope nevertheless to have shown that there are some remarkable parallells between sects like the ``Name-Worshippers' and the Early Modern ``enthoussiasts''. This again renders more clear why the destruction of these enthousisast tendencies throughout the European (Contra-)Reformation has had such marked out implications for mathematical invention and the place of imagination in philosophy and science in general. This line of research opens up some other interesting possibilities, like a re-assesment of the mathematics implied in the diagrammatic representations interspersed throughout the works of Giordano Bruno, or a re-evaluation of the philosophical and theological sources of Cantor's ideas on infinity, developed as they were in a Germany outliving the Romantic reaction against the excesses of Refomatory Enlightenment, as well as, say, Grothendieck's reflexions --- in the first part of {\em R\'{e}coltes et semailles} --- on ``Le Reveur'', his source of inspiration, during the period of the beat-generation and the ensuing psychedelc revolution in the West.

\end{document}